\documentclass{elsarticle}
\usepackage{amsmath}
\usepackage{amsfonts}
\usepackage{amssymb}
\renewcommand{\bf}{\bfseries}
\renewcommand{\it}{\itshape}
\def\a{\alpha}
\def\b{\beta}
\def\be{\begin{equation}}
\def\ee{\end{equation}}
\def\nn{\mathbb{N}}

\newtheorem{Theorem}{Theorem}[section]

\newtheorem{Corollary}[Theorem]{Corollary}

\def\proof#1. {\par
                      \ifdim\lastskip<15pt
                      \removelastskip\penalty-200
                      \vskip5pt plus3pt minus3pt
                      \fi
                       {\def\a{#1}
                       \ifx\a\empty
                       {\noindent\bf Proof.}
                       \else
                       {\noindent\bf Proof of #1.}
                       \fi}\enspace}

\def\endproof{\hfill\hspace{-6pt}\rule[-4pt]{6pt}{6pt}
\vskip8pt plus3pt minus 3pt}

\allowdisplaybreaks[4]
\journal{}
\begin{document}

\begin{frontmatter}

\title{Bounds for extreme zeros of some classical orthogonal polynomials}

\author{K Driver\fnref{label1}}
\ead{kathy.driver@uct.ac.za}\fntext[label1]{Research by this author is supported by
the National Research Foundation of South Africa under grant number
2053730.}
\address{Department of Mathematics and Applied Mathematics, University of Cape Town, Private Bag X3,
Rondebosch 7701, Cape Town, South Africa}
\author{K Jordaan \fnref{label2}}\ead{kjordaan@up.ac.za}\fntext[label2]{Research by this author is supported by the National Research Foundation of South Africa
under grant number 2054423.}
\address{Department of Mathematics and Applied Mathematics, University of Pretoria, Pretoria, 0002, South Africa}

\begin{abstract}
\noindent We derive upper bounds for the smallest zero and lower bounds for the largest zero of
Laguerre, Jacobi and Gegenbauer polynomials. Our approach uses mixed three term recurrence
relations satisfied by polynomials corresponding to different parameter(s) within the same classical
family. We prove that interlacing properties of the zeros impose restrictions on the possible
location of common zeros of the polynomials involved and deduce strict bounds for the extreme zeros
of polynomials belonging to each of these three classical families. We show numerically that the bounds
generated by our method improve known lower (upper) bounds for the largest (smallest) zeros of polynomials
in these families, notably in the case of
Jacobi and Gegenbauer polynomials.

\end{abstract}


\end{frontmatter}

\section{Introduction}

\noindent The behavior of the zeros of orthogonal polynomials has attracted significant
interest from both theoreticians and numerical analysts since the first results were proved
by Markov and Stieltjes in the late 19th century. A wide range of tools and techniques have been developed to analyse different properties of the zeros including Markov's theorem on the monotonicity of zeros (cf. \cite{Markov}); Sturm's comparison theorem for the zeros of solutions of second order differential equations (cf. \cite{Szego});  Obrechkoff's theorem on Descartes' rule of sign (cf. \cite{Obr}); and the Wall-Wetzel theorem on eigenvalues of Jacobi matrices (cf. \cite{Wall}). Discussion of results and techniques related to the properties of the zeros of both general orthogonal polynomials and classical families of orthogonal polynomials can be found in \cite{Elbert} and \cite{Muldoon}.

\medskip
\noindent It is well known that if $\{p_n\}_{n=0}^\infty$ is any orthogonal sequence, then the zeros of $p_n$ are real and simple and each open interval with endpoints at successive zeros of $p_n$ contains exactly one zero of $p_{n-1}$; a property called the interlacing of zeros. Stieltjes (cf. \cite{Szego}, Theorem 3.3.3) extended this interlacing property by proving  that if $m < n-1$, provided $p_m$ and $p_n$ have no common zeros, there
exist $m$ open intervals, with endpoints at successive zeros of
$p_n$, each of which contains exactly one zero of $p_m$. Beardon (cf. \cite{Beardon}, Theorem 5)
proved that one can say more, namely, for each $m < n-1$, if $p_m$ and $p_n$
are co-prime, there exists a real polynomial $S_{n-m-1}$ of
degree $n-m-1$ whose real simple zeros together with those of $p_m$, interlace with the zeros of $p_n$.
The polynomials $S_{n-m-1}$ are the dual polynomials introduced
by de Boor and Saff in \cite{DeB} or equivalently, the associated polynomials analysed by
Vinet and Zhedanov in \cite{Vin}.  Related results are proved by Segura in \cite{Segura}.

\medskip
\noindent In recent years, authors including Ismail and Muldoon (cf. \cite{Ism}), Krasikov (cf. \cite{Krasikov}), Gupta and Muldoon (cf. \cite{Gupta}), Ismail and Li (cf. \cite{ISLI}), and Dimitrov and Nikolov (cf. \cite{dimitrov}) have developed interesting methods including the use of chain sequences and the derivation of inequalities for real-root polynomials to refine and improve upper and lower bounds for extreme zeros of classical orthogonal polynomials. In this paper, we prove that if two polynomials $q_k$ and $p_{n+1}$, $k<n$, satisfy a three term recurrence relation of a certain type, common zeros of $q_k$ and $p_{n+1}$ cannot occur at consecutive zeros, nor at the largest or smallest zero, of $p_{n+1}.$  We show how this knowledge of the allowable location of common zeros, together with the Stieltjes interlacing property, leads to strict lower (upper) bounds for the largest (smallest) zero of $p_{n+1}$ in any orthogonal  sequence $\{p_n\}_{n=0}^\infty.$  We apply our results to three families of classical orthogonal sequences, namely, to sequences of Jacobi, Gegenbauer and Laguerre polynomials and compare our bounds to those derived by other authors for these classical families.

\medskip

\section{Results}

\noindent The first result that we use is proved in \cite{Jacobi}.

\begin{Theorem}\label{one} (cf. \cite[Lemma 3.2]{Jacobi})
Let $\{p_n\}_{n=0}^\infty$ be any sequence of polynomials
orthogonal on the (finite or infinite) interval $(c,d)$. Let
 $g_{n-1}$ be any polynomial of degree $n-1$ that satisfies
\begin{equation}\label{9}g_{n-1}(x)=a_n(x)p_{n+1}(x)-(x-A_n)b_n(x)p_n(x)\end{equation}
for for each fixed $n \in \nn$ where $A_n$ is a constant (depending on $n$) and $a_n$ and $b_n$ are functions with $b_n(x)\neq 0$ for $x\in(c,d)$.
Then, for each fixed $n\in\nn$,
\begin{itemize}
\item[(i)] the zeros of $g_{n-1}$ are all real and simple
and, together with the point $A_n$, they interlace with the zeros of
$p_{n+1}$ if $g_{n-1}$ and $p_{n+1}$ are co-prime;
\item[(ii)] if $g_{n-1}$ and $p_{n+1}$ are not co-prime, then

\noindent a) they have one common zero located at $x=A_n;$

\noindent b) the $n-1$ zeros of $g_{n-1}$ interlace with the $n$ remaining (non-common) zeros of $p_{n+1};$

\noindent c)  if $w_{n+1} < \dots < w_1$ denote the zeros of $p_{n+1},$ then $A_n = w_k$ for some $k \in \{2,\dots,n\}.$
\end{itemize}
\end{Theorem}

\begin{Corollary} \label{cor} For each fixed $n\in\nn,$ $A_n$ is a strict lower bound for the largest zero of $p_{n+1}$ and a strict upper bound for the smallest zero of $p_{n+1}$.
\end{Corollary}

\noindent {\bf{Remark}.} Equation (1) and the conditions of Theorem \ref{one} hold when $g_{n-1} =  p_{n-1}$ by virtue of the three term recurrence relation satisfied by any orthogonal sequence $\{p_n\}_{n=0}^\infty$ and hence Corollary \ref{cor} is a generalisation of Theorem 2 in \cite {BDJ} while Theorem 2.1 (ii) is a generalisation of  Theorem 3 in \cite {Beardon} since it accommodates the case when $p_{n-1}$ and $p_{n+1}$ have a common zero.

\medskip
\noindent Our main result is a generalisation of Theorem \ref{one}.

\begin{Theorem}\label{main}
Let $\{p_n\}_{n=0}^\infty$ be any sequence of polynomials
orthogonal on the (finite or infinite) interval $(c,d)$. Let
 $g_{n-k}$ be any polynomial of degree $n-k$ that satisfies, for each $n\in\nn,$ and $k \in \{1,\dots,n-1\},$
\begin{equation}\label{13}f(x)g_{n-k}(x)=G_k(x) p_n(x) + H(x) p_{n+1}(x)\end{equation}
where $f(x) \neq 0$ for $x\in(c,d)$ and $H(x)$, $G_k(x)$ are polynomials with deg($G_k)=k$.
Then, for each fixed $n\in\nn$,
\begin{itemize}
\item[(i)] the $n$ real, simple zeros of $G_k g_{n-k}$ interlace with the zeros of
$p_{n+1}$ if $g_{n-k}$ and $p_{n+1}$ are co-prime;
\item[(ii)] if $g_{n-k}$ and $p_{n+1}$ are not co-prime and have $r$ common zeros, counting multiplicity, then

\noindent a) $r \leq$ min $\{k,n-k\}$ ;

\noindent b) these $r$ common zeros are simple zeros of $G_k$;

\noindent c) no two successive zeros of $p_{n+1}$ , nor its largest or smallest zero, can also be zeros of $g_{n-k}$;

\noindent d) the $n-2r$ zeros of $G_k g_{n-k},$ none of which is also a zero of $p_{n+1}$, together with the $r$ common zeros of $g_{n-k}$ and $p_{n+1},$ interlace with the $n+1-r$ remaining (non-common) zeros of $p_{n+1}.$

\end{itemize}
\end{Theorem}

\proof{Theorem \ref{main}}. Let $w_{n+1} < \dots < w_1$
denote the zeros of  $p_{n+1}.$
\begin{itemize}
\item[(i)]

From (2), provided $p_{n+1}(x)\neq 0$, we
have \be\label{01} \frac{f(x)
g_{n-k}(x)}{p_{n+1}(x)} = H(x)
+\frac{G_{k}(x) p_{n}(x)}{ p_{n+1}(x)}.\ee

Further, \[\frac{p_n(x)}{p_{n+1}(x)}=\sum_{j=1}^{n+1}\frac{A_j}{x-w_j}\]
where  $A_j>0$ for  each $j\in\{1,\dots,n+1\}$ (cf. \cite[Theorem 3.3.5]{Szego}). Therefore (\ref{01}) can be written
as \be\label{04} \frac{f(x)
g_{n-k}(x)}{p_{n+1}(x)} = H(x)
+\sum_{j=1}^{n+1}\frac{G_{k}(x)A_j}{x-w_j},~~~ x\neq w_j.\ee Since
$p_{n+1}$ and $p_{n}$ are always
co-prime while $p_{n+1}$ and
$g_{n-k}$ are co-prime by assumption, it follows
from (2) that $G_k(w_j)\neq 0$ for  any
$j\in\{1,2,\dots,n+1\}.$ Suppose that $G_k$ does not change sign in the interval $I_j=(w_{j+1},w_j)$ where
$j\in\{1,2,\dots,n\}$. Since $A_j>0$ and the polynomial $H$ is bounded on $I_j$ while the right hand side of (\ref{04}) takes arbitrarily large positive and negative values on $I_j$, it follows that
$g_{n-k}$ must have an odd number of zeros in
each interval in which $G_k$ does not change sign. Since $G_k$ is of degree $k$, there are at least
$n-k$ intervals $(w_{j+1},w_j)$, $j\in\{1,\dots,n\}$ in which $G_k$ does not change sign and so each of these intervals must
contain exactly one of the $n-k$ real, simple zeros of $g_{n-k}$.
We deduce that the $k$ zeros of $G_k$ are real and simple and, together with the
zeros of $g_{n-k}$, interlace with the $n+1$ zeros of $p_{n+1}.$

\item[(ii)] If $r$ is the total number of common zeros of $p_{n+1}$ and
$g_{n-k}$ counting multiplicity then each of these $r$ zeros is a simple zero of $p_{n+1}$
and it follows from (2) that any common zeros of
$g_{n-k}$ and $p_{n+1}$ must also be zeros of $G_k$ since $p_n$ and
$p_{n+1}$  are co-prime. Therefore, $r\leq\mbox{min}\{k,n-k\}$ and there must be at least
$(n-2r)$ open intervals of the form $I_j=(w_{j+1},w_j)$ with endpoints at successive zeros of $p_{n+1}$
where neither $w_{j+1}$ nor $w_j$ is a zero of $g_{n-k}$ or $G_k(x)$.

If $G_k$ does not change sign in an interval $I_j=(w_{j+1},w_{j})$, it follows from (\ref{04}),
since $A_j>0$ and $H$ is bounded while the right hand side takes arbitrarily large
positive and negative values for $x\in I_j$, that $g_{n-k}$ must have an odd number of zeros
in that interval. Therefore, in at least $(n-2r)$ intervals $I_j$ either $g_{n-k}$ or $G_k$,
but not both, must have an odd number of zeros counting multiplicity. On the other hand, $g_{n-k}$ and $G_k$
have at most $(n-k-r)$ and $(k-r)$ real zeros respectively that are not zeros of $p_{n+1}$.
We deduce that there must be at most $(n-2r)$ intervals $I_j =(w_{j+1},w_j)$ with endpoints
at successive zeros $w_{j+1}$ and $w_j$ of $p_{n+1}$ neither of which is a zero of $g_{n-k}$.
It is straightforward to check that if the number of intervals $I_j =(w_{j+1},w_j)$
with endpoints at successive zeros of $p_{n+1}$ neither of which is a zero
of $g_{n-k}$ is equal to $n-2r$, this is only possible if no two consecutive zeros of $p_{n+1},$
nor the largest or smallest zero of $p_{n+1}$, can be common zero(s) of $p_{n+1}$ and $g_{n-k}$.
This proves a) to c) and d) follows from c).

\end{itemize}
\endproof

\begin{Corollary} \label{cor2}For each fixed $n\in\nn$ and each $k \in \{1,\dots,n-1\},$ the largest zero of $G_k$ is a strict lower bound for the largest zero of $p_{n+1}$. Similarly, the smallest zero of $G_k$ is a strict upper bound for the smallest zero of $p_{n+1}$.
\end{Corollary}

\noindent {\bf {Remarks}.}

\begin{itemize}

\item[1.] Equation (2) and the conditions of Theorem \ref{main} hold when $g_{n-k} =  p_{n-k}$. This can be seen (cf. \cite[Thm 4]{Beardon}) by iterating the three term recurrence relation satisfied by any orthogonal sequence $\{p_n\}_{n=0}^\infty$. Theorem \ref{main} (ii) is more general than Theorem 5 in \cite {Beardon} for orthogonal polynomials
since it caters for the case when $p_{n-k}$ and $p_{n+1}$ have common zeros.

\item[2.] Gibson (cf. \cite{Gibson}) proved that if $\{p_n\}_{n=0}^{\infty}$ is any orthogonal sequence and $m < n-1$, then $p_n$ and $p_m$ can have at most min$\{m,n-m-1\}$ common zeros and he showed that this upper bound is sharp. Theorem \ref{main} (ii) a) generalises Gibson's result.

\item[3.] It turns out that common zeros of two polynomials is a frequent occurrence within most classical orthogonal families (cf. \cite {Aharanov}) and also occurs when the two polynomials are members of different orthogonal sequences within the same family (cf. \cite{Geg}, \cite{Laguerre} and \cite{Jacobi}). It is therefore important to accommodate the possibility of common zeros in the analysis of Stieltjes interlacing of zeros of polynomials.
\end{itemize}

\section{Applications to classical families and comparison with existing bounds}

\noindent Markov's Theorem proves that the zeros of classical orthogonal polynomials vary monotonically with an increase (or decrease) in the parameter(s). Mixed three term recurrence relations of the type given in (\ref{9}) and (\ref{13}) provide functional equations linking polynomials from the same classical orthogonal family but corresponding to different values of the parameter(s). From the perspective of deriving good bounds for the extreme zeros, those mixed three term recurrence
relations that involve the largest possible parameter difference between the polynomials, or alternatively, no parameter difference but the largest possible difference in the degree, while still preserving interlacing properties for their zeros, are particularly useful.

\medskip
\noindent In the discussion below, we shall denote the zeros of the polynomial $p_{n+1}$ of degree $n+1$ by
$w_{n+1} <\dots < w_1$.

\begin{itemize}

\item [(i)] \begin{itemize} \item[(a)]The sequence of Laguerre polynomials $\{L_{n}^\alpha\} _{n=0}^\infty $ is orthogonal
on the interval $(0,\infty)$ with respect to the
weight function $x^\alpha e^{-x}$ for $\alpha > -1.$  Upper bounds for the smallest zero $w_{n+1}$
of the Laguerre polynomial $L_{n+1}^{\a}$ , $\a > -1$ , namely
\begin{equation}\label{GM0}w_{n+1}<\frac{(\a+1)(\a+2)}{n+\a+2}\end{equation} and
\begin{equation}\label{GM1}w_{n+1}<\frac{(\a+1)(\a+3)}{2n+\a+3}\end{equation}
were obtained by Hahn (cf. \cite{Hahn}) and Szeg\"{o} (cf. \cite[eqn. (6.31.12)]{Szego}) respectively. Gupta and Muldoon recovered these bounds in \cite[eqns. (2.9) and
(2.10)]{Gupta} and the bounds also follow from Corollary \ref{cor} and Theorem 1 (iv) and (v) in \cite{Laguerre}.

\medskip

Replacing $\a$ by $\a+1$, $n$ by $n-1$ in \cite[eqn. (13)]{Laguerre} we have
\begin{eqnarray*}x^4L_{n-2}^{\a+5}(x)&=&(\a+3)(n+\a+1)((\a+2)(\a+4)-(2n+\a+2)x)L_{n-1}^{\a+1}(x)\\&&+n((n-1)x^2+2(\a+3)(n-1)x-(\a+2)(\a+3)(\a+4))L_{n}^{\a+1}(x)\end{eqnarray*}
and using \cite[eqn. (4)]{Laguerre} together with
the three term recurrence relation for Laguerre polynomials (cf. \cite{Szego}) we obtain
\begin{eqnarray}\label{ex} x^5L_{n-2}^{\a+5}(x) &=& (n+ \a +1)\left((\a+1)_4 - (\a+2)_2 (3n + 2 \a+5) x + (n + \a + 2)_2 x^2\right) L_n^\a(x)\\&&\nonumber
-(n+1) \left((\a+1)_4 - (\a+2)_2 (3n + \a +1) x + (n-1) n x^2\right) L_{n + 1}^{\a}(x)\end{eqnarray} where \[(\alpha)_k= \alpha(\alpha+1)\ldots(\alpha+k-1),\quad k \geq1,\quad
k \in \nn\] is Pochhammer's symbol.
The smallest zero of the quadratic polynomial coefficient of $L_n^\a$ in (\ref{ex}) is
\begin{eqnarray} \frac{(\a+2)_2 (3n + 2 \a + 5) -
   \sqrt{(\a+2)_2(9 (\a+2)_2 +
       2 (2\a+5) (\a^2+5\a+10) n + (5\a^2+ 25 \a +38) n^2)}}{2 (n + \a + 2)_2}\label{ours} \end{eqnarray}
and this is a sharper upper bound for the smallest zero of $L_{n+1}^\a$ than the upper bound
\begin{equation}\label{GM2}w_{n+1}<\frac{( \a+1) (\a+2) (\a+4) (2n+\a+3)}{ (\a+1)^2 (\a+2) + (5 \a+11) ( n+1) (n+ \a + 2)}\end{equation} derived by Gupta and Muldoon in \cite[eqn. (2.11)]{Gupta} as illustrated in table \ref{lagu}.
\begin{table}[!h] \caption{Comparison of upper bounds for the smallest zero of $L_{n+1}^{\a}(x)$ for different values of $n$ and $\alpha$.}\label{lagu}
\begin{tabular}{|c|c|c|c|c|}
\hline
&&&&\\
Values of $\a$ and $n$& $\a=340.56$, $n=12$& $\a=65.3$, $n=20$& $\a=-0.9$, $n=100$ & $\a=-0.9$, $n=3$ \\[0.2cm]
\hline
& & & & \\
Smallest zero of $L_{n+1}^{\a}$&251.815653&27.67770854& 0.00103831& 0.0259078\\[0.2cm]
\hline
& & & & \\
New bound (\ref{ours})&275.856&32.6035&0.00103832&0.0259079\\[0.2cm]
\hline
& & & & \\
Bound (\ref{GM2}) given in \cite{Gupta}&309.594&36.6163&0.00103833&0.0259082\\[0.2cm]
\hline
& & & & \\
Bound (\ref{GM1}) given in \cite{Gupta}&319.258&41.8125&0.00103909&0.0259259\\[0.2cm]
\hline
& & & &  \\
Bound (\ref{GM0}) given in \cite{Gupta}&330&51.111&0.00108803&0.268293\\[0.2cm]
\hline
\end{tabular}\end{table}
 \newpage
\noindent{\bf{Remark.}} As observed in \cite {Gupta} , it is evident from (5), (6) and (9) that the upper bounds
for the smallest zero will be sharpest for $\alpha$ close to $-1$.

\medskip
\item[(b)]
The lower bound for the largest zero $w_{1}$ of $L_{n+1}^{\a}(x)$, $\a>-1$
\begin{eqnarray}
\label{Neumann}w_1&>&3n-1\end{eqnarray}
is given by Neumann in \cite{Neumann} while
\begin{eqnarray}\label{Bottema}
w_1&>&4n+4+\a-16\sqrt{2n+2}\end{eqnarray} in \cite{Bottema} is due to Bottema and
\begin{eqnarray}\label{L3}w_{1}&>&2n+\a+1\end{eqnarray}
was obtained by Szeg\H{o} (cf. \cite[eqn. (6.2.14)]{Szego}). Note that (\ref{L3}) also follows from Theorem 1 (i) in \cite{Laguerre} and Corollary \ref{cor}.

\medskip
By iterating the three term recurrence relation for Laguerre polynomials we obtain
\begin{eqnarray}\label{rec}\lefteqn{(\a + n-1) (\a + n)L_{n - 2}^{\a}(x)}\nonumber \\ &=& (x^2- 2 (2 n+\a) x +3n^2+3\a n + \a^2 + 1)L_{n}^{\a}(x) - (n+1) (2n+\a-1 - x)L_{n+1}^{\a}(x)\end{eqnarray} and the largest zero of the quadratic polynomial coefficient of $L_n^\a$ in (\ref{rec}) yields the
improved lower bound \begin{equation}\label{L4}w_{1}>2n+\a+\sqrt{n^2+\a n+1}\end{equation} for the largest zero of $L_{n+1}^{\a}$.

\noindent{\bf Remark.}
It is interesting to note that although (\ref{L4}) is an improvement on previous bounds for most values of $n$ and $\a$, the bound (\ref{Neumann}) due to Neumann, which does not depend on $\a$, compares favourably when both $n$ and $\a$ are small.

\begin{table}[!h] \caption{Comparison of lower bounds for the largest zero of $L_{n+1}^{\a}(x)$ for different values of $n$ and $\alpha$.}
\begin{tabular}{|c|c|c|c|c|}
\hline
&&&&\\
Values of $\a$ and $n$& $\a=340.56$, $n=12$& $\a=65.3$, $n=20$& $\a=-0.9$, $n=100$ & $\a=-0.9$, $n=4$  \\[0.2cm]
\hline
& & & & \\
Largest zero of $L_{n+1}^{\a}$&469.74362252&172.7701885&591.362&11.1263\\[0.2cm]
\hline
& & & & \\
New bound (\ref{L4})&429.612&146.616&298.654&10.7606\\[0.2cm]
\hline
& & & & \\
Bound (\ref{L3}) given in \cite{Szego}&365.56&106.3&200.1&8.1\\[0.2cm]
\hline
& & & &  \\
Bound (\ref{Bottema}) given in \cite{Bottema}&310.976&45.6081&175.69&-31.4964\\[0.2cm]
\hline
& & & &  \\
Bound (\ref{Neumann}) given in \cite{Neumann}&35&59&299&11\\[0.2cm]
\hline
\end{tabular}\end{table}
\newpage

\end{itemize}

\medskip
\item [(ii)]  The sequence of Jacobi polynomials
$\{P_{n}^{\alpha,\beta}\} _{n=0}^\infty $ is orthogonal on the interval $(-1,1)$
with respect to the weight function $(1-x)^\a(1+x)^\b$ for $\alpha$, $\beta>-1$.

It was proved in \cite[Thm 2.1(c)]{Jacobi} that for any fixed $n\in \nn$ and
$\alpha$, $\beta>-1$, (1) holds with

\medskip

\[g_{n-1}=P_{n-1}^{\alpha+4,\beta},~
A_n = \frac{2n(n+\a+\b+3)+(\a+3)(\beta-\a)}
{2n(n+\a+\b+3)+(\a+3)(\a+\beta+2)}~\mbox{ and }~p_{n+1} = P_{n+1}^{\alpha,\beta}.\]
It follows from Corollary \ref{cor} that

\begin{eqnarray}\label{J1}w_1& > & 1 - \frac{2(\a+1)(\a+3)}
{2n(n+\a+\b+3)+(\a+3)(\a+\beta+2)} \end{eqnarray}

\medskip

Similarly, from \cite[Cor 2.2(c)]{Jacobi}, (1) holds with

\medskip
\[g_{n-1}=P_{n-1}^{\alpha,\beta+4},~A_n=-\frac{2n(n+\a+\b+3)+(\b+3)(\a-\b)}
{2n(n+\a+\b+3)+(\b+3)(\a+\beta+2)}~\mbox{ and }~p_{n+1}=P_{n+1}^{\alpha,\beta}.\]

\medskip

and we deduce that

\begin{eqnarray}\label{J2}w_{n+1}& <& -1 + \frac{2(\b+1)(\b+3)}
{2n(n+\a+\b+3)+(\b+3)(\a+\beta+2)} \end{eqnarray}

\medskip

\noindent{\bf Remark.} It should be noted that for any fixed $n\in \nn$, (\ref{J1}) is a particularly
sharp lower bound for the largest zero of $P_{n+1}^{\a,\b}$ when $\a$ is close to $-1$ and similarly
(\ref{J2}) provides a sharp upper bound for the smallest zero of $P_{n+1}^{\a,\b}$ when $\b$ is
close to $-1$. Asymptotically in $n$, (\ref{J1}) and (\ref{J2}) yield $w_1 > 1- O(\frac 1{n^2})$ and
$w_{n+1} < - 1 + O(\frac 1{n^2})$ respectively.

\medskip
In table \ref{Jacobi} we compare (\ref{J1}) to the bounds given by Szeg\"{o} in
\cite[eqns. (6.2.11) and (6.2.12)]{Szego}
\begin{eqnarray}\label{z1}w_{1}&>&\frac{2n+\beta-\a}{2n +\a+\b+2}\\
\label{z2}w_1&>&\frac{n}{n+\a+1}~\mbox{ for }\b\geq\a.\end{eqnarray}

\begin{table}[!h]\label{Jacobi} \caption{Comparison of lower bounds for the largest zero of $P_{n+1}^{\a,\b}(x)$ for different values of $n$, $\alpha$ and $\b$.}
\begin{tabular}{|c|c|c|c|c|c|c|}
\hline
&&&&&&\\
Value of $n$&$n=3$&$n=15$&$n=11$&$n=18$&$n=3$&$n=3$\\[0.2cm]
\hline
&&&&&&\\
Value of $\a$ & $\a=-0.9$& $\a=-0.9$& $\a=30.9$& $\a=30.9$&$\a=30.9$&$\a=0.9$  \\[0.2cm]
\hline
&&&&&&\\
Value of $\b$& $\b=-0.8$& $\b=-0.8$&$\b=-0.8$ &$b=32.8$& $\b=32.8$& $\b=132.8$\\[0.2cm]
\hline
& & &&&& \\
Largest zero of $P_{n+1}^{\a,\b}$&0.984119&0.999143&0.141417&0.71025&0.30295&0.989891\\[0.2cm]
\hline
& & &&&&  \\
New bound (\ref{J1})&0.984109&0.999142&-0.0507338&0.590098 &0.182432&0.989162\\[0.2cm]
\hline
& & & &&&\\
Bound (\ref{z1}) given in \cite{Szego}&0.968254&0.993399&-0.179298& 0.372665&0.110181&0.973183\\[0.2cm]
\hline
&&&&&&\\
Bound (\ref{z2}) given in \cite{Szego}&0.967742&0.993377&NA&0.360721&0.0859599& 0.612245\\
\hline
\end{tabular}\end{table}
\newpage

\item [(iii)] The sequence of Gegenbauer (or ultraspherical) polynomials
$\{C_{n}^{\lambda}\}_{n=0}^\infty $ is the one-parameter classical family obtained from the
sequence of Jacobi polynomials $\{P_{n}^{\alpha,\beta}\} _{n=0}^\infty $
with $\alpha = \beta = \lambda - \frac12$. The sequence is orthogonal on $[-1,1]$ with respect to
the weight function $(1-x^2)^{\lambda-\frac 12}$ for $\lambda>-\frac 12$ ,
$n \in {\mathbb N}$. The zeros of $C_n^\lambda$ lie in the open
interval $(-1,1)$ and are symmetric about the origin with a simple
zero at the origin when $n$ is odd. From \cite[Theorem 2 (ii) d)]{Geg} and Corollary \ref{cor2}, it follows that

\begin{eqnarray}\label{geg1}w_1^2 &>& 1 - \frac{(2\lambda+1)(2\lambda+3)}
{ n(n+2\lambda+2)+(2\lambda+1)(2\lambda+3)}.\end{eqnarray}

From (\ref{J1}), letting $\a=\b=\lambda-\frac{1}{2}$, we obtain
\begin{equation}\label{geg0}w_1 > 1 - \frac{(2\lambda+1)(2\lambda+5)}
{4 n(n+2\lambda+2)+(2\lambda+1)(2\lambda+5)}\end{equation}

\medskip

\noindent{\bf Remark.} The bounds obtained in (\ref{geg1}) and (\ref{geg0}) are particularly sharp when
$\lambda$ is close to $-\frac 12$.  For $n$ large, (\ref{geg1}) and (\ref{geg0})
yield $w_1^2 > 1- O(\frac 1{n^2})$ and $w_1 > 1- O(\frac 1{n^2})$ respectively.

\medskip

In table 4 we compare (\ref{geg1}) and (\ref{geg0}) to the bound given by Szego (cf. \cite[eqn. (6.2.13)]{Szego})
\begin{equation}\label{geg2}w_1^2\geq {1-\frac{2\lambda+1}{n+2\lambda+1},}\end{equation} noting that
this lower bound also follows immediately from \cite[Theorem 2 (ii) c)]{Geg} by virtue of Corollary \ref{cor2}.
\end{itemize}
\begin{table}[!h]\label{gegg} \caption{Comparison of lower bounds for the largest zero of $C_{n+1}^{\lambda}(x)$ for different values of $n$ and $\lambda$.}
\begin{tabular}{|c|c|c|c|c|}
\hline
&&&&\\
Value of $n$ and $\lambda$&$n=30$, $\lambda=-0.49$&$n=3$, $\lambda=-0.4956$&$n=50$, $\lambda=30.9$&$n=5$, $\lambda=60.49$\\[0.2cm]
\hline
& & && \\
Largest zero of $C_{n+1}^{\lambda}$& 0.999978&0.999267&0.897416&0.289296\\[0.2cm]
\hline
& & && \\
New bound (\ref{geg0})& 0.999978&0.999267&0.844369&0.142782\\[0.2cm]
\hline
& & && \\
New bound (\ref{geg1})& 0.999978&0.999266&0.763561&0.201482\\[0.2cm]
\hline
& & & &\\
Bound (\ref{geg2}) given in \cite{Szego}&0.999667&0.998537&0.66578&0.198435\\[0.2cm]
\hline
\end{tabular}\end{table}
\newpage

\section{Concluding remarks}

\begin{itemize}

\item[1.] The bounds derived in this paper improve previously known upper (lower) bounds for the smallest (largest) zeros of Laguerre, Jacobi and Gegenbauer polynomials, significantly so in the cases of both Jacobi and Gegenbauer. Our results complement those obtained by other authors who derive upper (lower) bounds for the largest (smallest) zeros, notably, van Doorn in \cite{vanDoorn}; Dimitrov and Nikolov in \cite{dimitrov} for the zeros of Laguerre, Jacobi and Gegenbauer polynomials; Ismail and Li in \cite{ISLI}; and {Krasikov} who gives uniform bounds for the extreme zeros of Laguerre polynomials in \cite{Krasikov}.

\item[2.] Various versions of (2) can be used to obtain sharper upper (lower) bounds for the smallest (largest) zeros of Laguerre, Jacobi and Gegenbauer polynomials by letting $k = 3,4,...$ The calculations become more complicated as the degree of the coefficient polynomial $G_k$ increases.
\end{itemize}

\end{document}